\documentclass[12pt]{article}
\usepackage{xcolor,graphicx} 
\usepackage{amsmath,amsthm,amssymb}

\usepackage{ifthen}

\usepackage{amsmath,amsthm,amssymb,mathtools}
\usepackage{pgf}
\usepackage{graphicx}
\usepackage{tikz}
\usetikzlibrary{backgrounds, patterns, shapes, calc, decorations.pathreplacing, decorations.pathmorphing}

\tikzset{%
    vtx/.style={draw,circle,very thin,white,fill=black,inner sep=0pt,text width=8pt},
    pth/.style={draw,decorate,decoration={snake},thick},
    edg/.style={ultra thick},
    smedg/.style={very thin},
    polyg/.style={very thin,fill=white},
    dotedg/.style={ultra thick,dotted},
    ellip/.style={very thin,black,fill=white},
    spedg/.style={rounded corners,line width=5pt,purple!50}, 
}
\tikzset{/pgf/foreach/parse=true}
\pgfdeclarelayer{backback}
\pgfdeclarelayer{back}
\pgfdeclarelayer{fore}
\pgfdeclarelayer{forefore}
\pgfsetlayers{backback,back,main,fore,forefore}

\newtheorem{theorem}{Theorem}

\newtheorem{lemma}[theorem]{Lemma}
\newtheorem{corollary}[theorem]{Corollary}

\newtheorem{example}[theorem]{Example}
\newtheorem{conjecture}{Conjecture}

\newtheorem{claim}[theorem]{Claim}

\newtheorem{obs}[theorem]{Observation}
\def\beq{\begin{equation}}\def\eeq{\end{equation}}
\def\qed{\ifhmode\unskip\nobreak\fi\quad\ifmmode\Box\else$\Box$\fi}

\def\qed{\hfill $\Box$}

\title{2-reachable subsets in two-colored graphs}
\author{Andr\'as Gy\'arf\'as\thanks{Alfr\'ed R\'enyi Institute of Mathematics, Budapest, P.O. Box 127, Budapest, Hungary, H-1364.
\texttt{gyarfas.andras@renyi.hu}, \texttt{sarkozy.gabor@renyi.hu}} \thanks{Research supported in part by
NKFIH Grant No. K132696.} \and G\'{a}bor N. S\'ark\"ozy\footnotemark[1]
\thanks{Computer Science Department, Worcester Polytechnic Institute, Worcester, MA.} \thanks{Research supported in part by
NKFIH Grants No. K132696, K117879.}}
\date{}

\begin{document}
\maketitle
\begin{abstract} A subset $X$ of vertices in a graph $G$ is a {\em diameter 2 subset} if the distance of any two vertices of $X$ is at most two {\em in $G[X]$}.  Relaxing this notion, a subset $X$ of vertices in a graph $G$ is a {\em 2-reachable subset} if the distance of any two vertices of $X$ is at most two {\em in $G$}. Related to recent attempts to strengthen a well-known conjecture of Ryser, English et al. conjectured that the vertices of a $2$-edge-colored cocktail party graph (the graph obtained from a complete graph with an even number of vertices by deleting a perfect matching) can be covered  by the vertices of two monochromatic  diameter $2$ subsets. In this note we prove the relaxed form of this conjecture, replacing diameter $2$ by $2$-reachable.  An immediate corollary is that $2$-colored cocktail party graphs on $n$ vertices must contain a monochromatic $2$-reachable subset with at least $n\over 2$ vertices (and this is best possible).
\end{abstract}

\section{Introduction}

In a $2$-edge colored graph $G=(V,E)$, let the two colors be color 1 and 2, i.e. $E=E_1\cup E_2, E_1\cap E_2=\emptyset$, and let $G_1=(V,E_1), G_2=(V,E_2)$.  In particular, $S_1(v), S_2(v)$ denote stars with center $v$ in colors 1 and 2, respectively.  Similarly, $N_1(v), N_2(v)$ denote the set of vertices adjacent to $v$ in colors 1 and 2, respectively.  The color of an edge $(x,y)$ is denoted by $col(x,y)$. Sometimes for convenience we think of color 1 as red and color 2 as blue.

A subset $X$ of vertices in a graph $G$ is a {\em diameter 2 subset} if the distance of any two vertices of $X$ is at most two {\em in $G[X]$}.  Relaxing this notion, a subset $X$ of vertices in a graph $G$ is called a {\em 2-reachable subset} if the distance of any two vertices of $X$ is at most two {\em in $G$} (thus we might use a middle vertex outside of $X$). This relaxation originated from a Ramsey-type question of Erd\H os and Fowler: is it true that any $2$-colored complete graph $K_n$ has a monochromatic diameter $2$ subgraph with at least $\lceil {3n\over 4} \rceil$ vertices? The relaxed question (with a short proof) was solved in \cite{GYfruit}. Using  a more elaborate form of that proof method, Erd\H os and Fowler \cite{EF} answered affirmatively the original question as well.

The covering version of the Ramsey-type question above is motivated by a possible strengthening of a classical conjecture attributed to Ryser (stated in the thesis of Henderson \cite{HE}).
The version below is from \cite{GYsztaki} (for details, see \cite{DEB}, page 3).


\begin{conjecture}\label{HR}
Let $r\geq 2$ be an integer and $G=(V,E)$ a graph with $\alpha(G)=\alpha$. Then in every $r$-edge coloring of $G$ there is a cover of $V(G)$ by $(r-1)\alpha$ monochromatic connected components.
\end{conjecture}

Mili\'{c}evi\'{c} in \cite{MI2} then DeBiasio et al. in \cite{DEB} conjectured that the diameter of the monochromatic components in Conjecture \ref{HR} can be bounded (for some details and partial results see Section \ref{conclusion}). The best bound for the diameter (namely 3) for the $\alpha=1,r=2$ case comes from the folklore result that any $2$-colored complete graph has a spanning monochromatic subgraph of diameter at most three.  For graphs with $\alpha=2,r=2$  we need a cover by two monochromatic graphs. If there is no missing edge at a vertex, the two monochromatic stars provide a cover by two diameter $2$ sets. The first graph when this trivial solution does not work is the ``cocktail party'' graph $G^c=(V,E)$, obtained by deleting a perfect matching from a complete graph with an even number of vertices. For any vertex $v\in V$ the unique vertex non-adjacent to $v$ in $G^c$ is denoted by $v'$. (Note that $G^c$ can also be viewed as a complete partite graph, where each partite class has size 2.) Surprisingly, it does not seem easy to decide whether a cover by two monochromatic diameter $2$ sets exists. English, McCourt, Mattes and Phillips conjectured a positive answer.

\begin{conjecture}\label{2diamcov}(\cite{EMMP}).
 In every $2$-coloring of the edges of the cocktail party graph $G^c=(V,E)$, there exist $A, B\subseteq V$ and two colors $i, j\in[2]$ such that $A\cup B=V$ and $A,B$ are diameter 2 subsets of $G^c_i,G^c_j$, respectively.
\end{conjecture}

In this note we prove the relaxed version of Conjecture \ref{2diamcov}.

\begin{theorem}\label{2reachcov}  In every $2$-coloring of the edges of the cocktail party graph $G^c=(V,E)$, there exist $A, B\subseteq V$ and $i, j\in[2]$ such that $A\cup B=V$ and $A,B$ are 2-reachable subsets of $G^c_i,G^c_j$, respectively.
\end{theorem}

It is natural to ask what happens with the Ramsey-type problem mentioned above if $2$-colored complete graphs are replaced by $2$-colored cocktail party graphs. Theorem \ref{2reachcov} has the following immediate corollary.

\begin{corollary}\label{cor} Any $2$-colored $G^c$ with $n$ vertices contains a monochromatic  $2$-reachable subset with at least  $n\over 2$ vertices.
\end{corollary}

The bound $n\over 2$ of Corollary \ref{cor} seems weak compared to ${3n\over 4}$ in case of complete graphs.  However, it is best possible, as the following simple example shows.

\begin{example}\label{exsharp} Let $G^c=(V,E)$, where $|V|=n$.  Write $V=X\cup Y$ where $X,Y$ are complete subgraphs in $G^c$. Color $G^c[X],G^c[Y]$ red and all other edges blue.
\end{example}

In this example there is no monochromatic $2$-reachable (let alone diameter $2$) $S\subset V$ with $|S|>{n\over 2}$.  Indeed, otherwise by the pigeonhole principle there are nonadjacent vertices $v,v'\in S$ and their distance is larger than two in both colors (infinite in red and $3$ in blue), a contradiction.

We use the following notation.  Let ${\cal{C}}_5^+$ be the family of graphs obtained from a five-cycle where vertices are replaced with independent sets and the edges between the parts corresponding to edges of the five-cycle are replaced by complete bipartite graphs. Note that the diameter of any graph in ${\cal{C}}_5^+$ is two.

\section{Proof of Theorem \ref{2reachcov}}

From now on $G^c=(V,E)$ is a $2$-colored cocktail party graph. A pair of vertices $x,y\in A\subseteq V$ (not necessarily an edge) is {\em critical for color $i\in [2]$} if the distance of $x,y$ in $G^c_i$ is greater than two.
Note that when $A=V$, a critical pair for color $i$ witnesses that $G^c_i$ has diameter greater than two. Thus we must have at least one critical pair in $V$ for both colors, otherwise $V$ itself has diameter at most two in one of the colors, satisfying the requirement of Theorem \ref{2reachcov}. Let $e=(x,y), f=(p,q)$ be critical pairs for colors 1 and 2, respectively.

\begin{obs}\label{easyobs} If   $e=(x,y)\notin E$  or $f=(p,q)\notin E$, then $V$ can be covered by two monochromatic stars of the same color. If  $e=(x,y)\in E$, then  $e\in E(G^c_2)$. If $f=(p,q)\in E$, then  $f\in E(G^c_1)$.
\end{obs}
Indeed,  the criticality of  $e=(x,y)\notin E$ means that there is no color 1 path $x,z,y$ in $G^c$ thus $S_2(x)\cup S_2(y)=V$. The argument is similar for $f=(p,q)$. The statements for $e,f\in E$ are immediate from the definition of critical pairs.
Thus we may always assume that $e\in E(G^c_2)$ and $f\in E(G^c_1)$, otherwise we are done.

The next lemma is the backbone of the proof; we show that we may assume that $e$ and $f$ are disjoint.

\begin{lemma}\label{disjoint}   If  $e\cap f\ne \emptyset$, then $V$ can be covered by two monochromatic stars or by two monochromatic graphs from ${\cal{C}}_5^+$.
\end{lemma}

Before proving Lemma \ref{disjoint} we show how to finish the proof of Theorem \ref{2reachcov} from it. Consider all critical edges of $G^c$, let $P\subseteq V, Q\subseteq V$ denote the vertex sets covered by the critical edges of colors 1 and 2, respectively. Then we apply  Lemma \ref{disjoint} for all pairs of critical edges such that one is inside $P$, the other is inside $Q$. We have two possibilities. If Lemma \ref{disjoint} provides the required covering of $G^c$ (by stars or ${\cal{C}}_5^+$-s) for {\em any} pair, the proof is finished. Otherwise $P\cap Q=\emptyset$, thus for $A=V\setminus P, B=V\setminus Q$ we have $A\cup B=V$. Now we claim that $A,B$ provide the required covering of $V$. Indeed, there is no critical edge of color 1 inside $A$, thus $A$ is 2-reachable in color 1. Similarly $B$ is 2-reachable in color 2, finishing the proof. Therefore it is enough to prove Lemma \ref{disjoint}.
\smallskip

\noindent{\bf Proof of Lemma \ref{disjoint}:} Suppose w.l.o.g.  that $e\cap f=\{x\}$, say $x=p$. (From Observation  \ref{easyobs} we know that $|e\cap f|=1$).
Define $X,Y\subset V$ as (see Figure 1)
$$  X=N_2(x), Y=N_2(y)\setminus (N_2(x)\cup \{x'\}).$$
Note, that since $e$ is critical in color $1$, we have $X\cap Y=\emptyset$ and
\beq\label{partition}
V= \{x, x', y, y'\} \cup X \cup Y.
\eeq

We start with the following claim.
\begin{claim} $q=y'$.
\end{claim}

Indeed,  $(y,q)\in E(G^c)$ would contradict the definition of $e,f$, since if $(y,q)$ has color 1, then $e$ is not critical in color 1 and if 
$(y,q)$ has color 2, then $f$ is not critical in color 2.
Thus $(y,q)\notin E(G^c)$, implying $q=y'$.


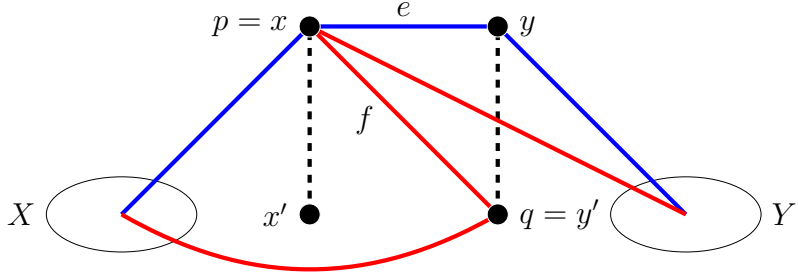
\begin{figure}
		       \centering
                \begin{tikzpicture}[scale=2.5]
                    \centering
				\useasboundingbox (0,-0.2) rectangle (1,1.5);
				\coordinate (a) at (0,0);
				\coordinate (b) at (1,0);		
				\coordinate (c) at (0,1);
				\coordinate (d) at (1,1);
				\coordinate (e) at (2,0);
			    \coordinate (f) at (-1,0);
    \def\bigrad{0.4};
    \def\smrad{0.2};
         \draw[ellip] (e) ellipse ({\bigrad} and {\smrad});
         \draw[ellip] (f) ellipse ({\bigrad} and {\smrad});

						\begin{pgfonlayer}{forefore}
					\draw (a) node[vtx,label={left:$x'$}] {};
					\draw (b) node[vtx,label={right:$q=y'$}] {};
					\draw (c) node[vtx,label={left:$p=x$}] {};
					\draw (d) node[vtx,label={right:$y$}] {};
                    \node [right] at (2.4,0) {$Y$};
                    \node [left] at (-1.4,0) {$X$};
                    \node [above] at (0.5,1) {$e$};
                    \node [left] at (0.4,0.5) {$f$};

  					\end{pgfonlayer}
				\begin{pgfonlayer}{fore}
					\draw[edg,dashed] (a) -- (c);
					\draw[edg,red] (c) -- (b);
					\draw[edg,blue] (c) -- (d);
					\draw[edg,dashed] (b) -- (d);
		            \draw[edg,blue] (d) -- (e);
                    \draw[edg,blue] (c) -- (f);
                    \draw[edg,red] (c) -- (e);
                    \path (-1,0) edge[red,bend right,ultra thick] (1,0);
							\end{pgfonlayer}
			\end{tikzpicture}
			\caption{The structure in the proof of Lemma \ref{disjoint}, where color 1 is red, color 2 is blue}
            \label{P_4^2}
	\end{figure}

From the definition of $Y$, $N_2(x)\cap Y=\emptyset$, so we get
\begin{obs} $Y\cup \{y'\}=N_1(x)$.
\end{obs}
Also, from the criticality of $f$ we have no path of length two from $p$ to $q=y'$ in color $2$, so we get
\begin{obs} $X\subseteq N_1(y')$.
\end{obs}

The structure of $G^c$ is shown in Figure 1.
Consider $x'$ and note that $x'$ is connected to all vertices in $X\cup Y$.
We partition $X,Y$ as follows. For $i\in [2]$  set
$$X_i=X\cap N_i(x'),Y_i=Y\cap N_i(x').$$

Next we start looking for coverings by two monochromatic stars. It is easy to check the following statements.

\begin{itemize}
\item (i)  if  $col(x',y)=2$ then $S_1(y')\cup S_2(y)=V,$
\item (ii) if  $col(x',y')=1$ then $S_1(y')\cup S_2(y)=V.$
\end{itemize}

From (i),(ii) we may assume that $col(x',y)=1,col(x',y')=2$ (note that these must be edges). Then, using this, we look at the effect of a potential degeneracy in the partitions $X_1,X_2$ and $Y_1,Y_2$.

\begin{itemize}
\item (iii)   if $X_1=\emptyset$ then $V=S_2(x')\cup S_2(y),$
\item (iv)  if $Y_1=\emptyset$ then $V=S_2(x')\cup S_2(x),$
\item (v) if $X_2=\emptyset$ then $V=S_1(x')\cup S_1(x),$
\item (vi)  if $Y_2=\emptyset$ then $V=S_1(x')\cup S_1(y').$
\end{itemize}
From (iii)-(vi) we may assume that the partitions $X_1,X_2$ and $Y_1,Y_2$ are not degenerate. From these assumptions the coloring of $G^c$ is restricted enough to conclude the proof with the following covering of
$V$ by two  monochromatic ${\cal{C}}_5^+$-s (note that this is the only place in the proof of Lemma \ref{disjoint} where we are using ${\cal{C}}_5^+$-s instead of stars):
$$G^c_2[x,y,Y_2,x',X_2,x] \mbox{ and }G^c_1[x,y',X_1,x',Y_1,x],$$
finishing the proof of the lemma (and Theorem \ref{2reachcov}).  \qed

\section{Concluding remarks}\label{conclusion}

For $r = 2$ Conjecture \ref{HR} is equivalent to K\"onig’s theorem \cite{KO} and the $r = 3$ case was proved by Aharoni \cite{AH}. For $\alpha(G) = 1$ (i.e. for complete graphs) the conjecture holds for $r=4$ (\cite{GYsztaki}, \cite{TU1}) and for $r=5$ (\cite{TU2}).

Mili\'{c}evi\'{c} \cite{MI2} conjectured that for $\alpha=1$ there exists a smallest $f=f(r)$ bounding the diameter of the monochromatic components in Conjecture \ref{HR}. It is known that $f(2)=3$ (folklore), the bound $f(3)\le 8$ in \cite{MI1} is improved to  $f(3)\le 4$ in \cite{DEB}.  This conjecture is extended for every $\alpha$ in \cite{DEB} defining $f=f(r,\alpha)$ as the smallest bound for the diameter of the monochromatic components in Conjecture \ref{HR}. It is known that $f(2,\alpha)\le \alpha^2+12\alpha+4$ \cite{DGHS}. The bound $f(2,2)\le 6$ in \cite{DEB} is improved to 4 in \cite{GS} but it is possible that the best bound is 3.

For certain graphs with $\alpha=2$, including odd antiholes, there is a covering by two monochromatic subgraphs of diameter at most three in every $2$-colorings \cite{GS}. However, it does not seem easy to decide whether the same is true for graphs defined by vertex disjoint odd antiholes with all possible edges between them.


\begin{thebibliography}{99}

\bibitem{AH} R. Aharoni. Ryser’s conjecture for tripartite 3-graphs, {\em Combinatorica}, {\bf 21}(1) (2001), pp. 1-4.

\bibitem {DGHS} L. DeBiasio, A. Girao, P. Haxell, M. Stein, A bounded diameter strengthening of K\"onig's theorem, to appear in {\em SIAM Journal on Discrete Mathematics}, arXiv:2409.18250

\bibitem{DEB} L. DeBiasio, Y. Kamel, G. McCourt, H. Sheats, Generalizations and strengthenings of Ryser's conjecture, {\em Electronic Journal of Combinatorics}, {\bf 28} (2021), P4.37.

\bibitem{EF} P. Erd\H os and T. Fowler, Finding large $p$-colored diameter two subgraphs, {\em Graphs and Combinatorics}, {\bf 15} (1999), pp. 21-27.


\bibitem{EMMP} S. English, G. McCourt, C. Mattes, M. Phillips, Low diameter monochromatic covers of complete bipartite graphs, {\em Journal of Combinatorics}, {\bf 15}(2) (2024), pp. 139-157.

\bibitem{HE} J. R. Henderson, Permutation decomposition of $(0-1)$-matrices and decomposition transversals, Ph.D. Thesis, California Institute of Technology, 1971, https://core.ac.uk/download/pdf/11814656.pdf.


\bibitem{GYsztaki} A. Gy\'arf\'as, Partition coverings and blocking sets in hypergraphs, {\em Communications of the Computer and Automation Institute of the Hungarian Academy of Sciences} {\bf 71} (1977), in Hungarian.

\bibitem{GYfruit} A. Gy\'arf\'as, Fruit salad, {\em Electronic Journal of Combinatorics}, {\bf 4} (1997), R8.

\bibitem{GS} A. Gy\'arf\'as, G. N. S\'ark\"ozy, Bounded diameter variations of Ryser's conjecture, arXiv:2505.02564

\bibitem{KO} D. K\"onig, Graphs and matrices, {\em Matematikai \'es Fizikai Lapok}, {\bf 38} (1931), pp. 116–119.



\bibitem{MI1} L. Mili\v cevi\v c, Commuting contractive families, {\em Fundamenta Mathematicae}, {\bf 231} (2015), pp. 225-272.

\bibitem{MI2} L. Mili\v cevi\v c, Covering complete graphs by monochromatically bounded sets, {\em Applicable Analysis and Discrete Mathematics}, {\bf 13} (2019), pp. 85-110.
\bibitem{TU1} Zs. Tuza, Some special cases of Ryser’s Conjecture, unpublished manuscript, 1979.
\bibitem{TU2} Zs. Tuza, Ryser’s conjecture on transversals of r-partite hypergraphs, {\em Ars Combinatoria}, {\bf 16} (1983), pp. 201-209.

\end{thebibliography}
\end{document}